\documentclass[secthm,final,amsthm]{elsart}
\usepackage{amssymb,amsmath}
\usepackage[all,2cell,ps]{xy}

\newcommand{\edge}[1]{\ar@{-}[#1]}
\newcommand{\node}{*+[o][F-]{ }}

\newcommand{\F}{\mathbb{F}}
\newcommand{\N}{\mathbb{N}}
\newcommand{\id}{\operatorname{id}}
\newcommand{\mc}{\mathcal}
\newcommand{\End}{\operatorname{End}}
\newcommand{\Hom}{\operatorname{Hom}}
\newcommand{\Mat}{\operatorname{Mat}}
\newcommand{\mcSR}{\mathcal{SR}}
\newcommand{\modsim}{/\hspace{-1mm}\sim\,}
\renewcommand{\phi}{\varphi}
\renewcommand{\theta}{\vartheta}

\begin{document}
\begin{frontmatter}
  
  \title{Classification of finite congruence-simple semirings with
    zero\thanksref{snf}}
 
  \author{Jens Zumbr\"agel}
  
  \thanks[snf]{This work has been supported by the Swiss National
    Science Foundation under grant no.  107887.}
  
  \ead{jzumbr@math.unizh.ch}
  
  \address{Institut f\"ur Mathematik, Universit\"at Z\"urich\\
    Winterthurerstrasse 190, 8057 Z\"urich, Switzerland}
  
  \begin{abstract}
    Our main result states that a finite semiring of order $>2$ with
    zero which is not a ring is congruence-simple if and only if it is
    isomorphic to a `dense' subsemiring of the endomorphism semiring
    of a finite idempotent commutative monoid.
    
    We also investigate those subsemirings further, addressing e.g. the
    question of isomorphy.
  \end{abstract}

\end{frontmatter}
  
\section{Introduction and main result}
  
Semirings, introduced by Vandiver~\cite{va34} in 1934, generalize the
notion of noncommutative rings in the sense that negative elements do
not have to exist. Since then there has been an active area of
research in semirings, both on the theoretical side and on the side of
applications e.g. in theoretical computer science. The reader may
consult the monographs of Golan~\cite{go99} and
Hebisch/Weinert~\cite{he93a} for a more elaborate introduction to
semirings.
 
In order to develop a structure theory for semirings, special interest
lies in semirings which are congruence-simple, meaning simple in the
sense that there are only trivial quotient semirings (see below for
precise definitions). The classification of simple commutative
semirings was achieved only recently in~\cite{el01}. In the general
case it has been shown later~\cite{mo04} that any finite simple
semiring of order $>2$ which is not a ring has to have either trivial
or idempotent addition.  In this paper we give a full classification
of finite simple semirings assuming they have a zero element.

\begin{defn}
  A set $R$ with two binary operations $+$ and $\cdot$ is called a
  \emph{semiring} (with zero) if $(R,+)$ is a commutative monoid,
  $(R,\cdot)$ is a semigroup, and the distributive laws $\ 
  x\cdot(y+z)=x\cdot y+x\cdot z\ $ and $\ (x+y)\cdot z=x\cdot z+y\cdot
  z\ $ hold for all $x,y,z\in R$; furthermore, the neutral element $0$
  of $(R,+)$ has to satisfy $0\cdot x=x\cdot 0=0$ for all $x\in R$,
  and is called \emph{zero}.
\end{defn}
  
We sometimes write the multiplication as concatenation, i.e.
$xy:=x\cdot y$ for $x,y\in R$. If $(R,\cdot)$ has a neutral element
$1\in R$ we call it \emph{one}. By a \emph{subsemiring} of a semiring
$R$ we mean a subset $S\subseteq R$ with $0\in S$ which is closed
under addition and multiplication. Naturally, $S$ itself is a
semiring.
  
As mentioned above our notion of simplicity relies on congruences.
  
\begin{defn}
  An equivalence relation $\sim$ on a semiring $R$ is called
  \emph{congruence} if \[x\sim y\quad\text{implies}\quad a+x\sim a+y,\ 
  ax\sim ay,\ xa\sim ya,\quad\text{for all }x,y,a\in R.\] The semiring
  $R$ is called \emph{congruence-simple} if its only congruences are
  $\sim\;=\id_R$ and $\sim\;=R\times R$.
\end{defn}

\begin{rem}
  Given a congruence $\sim$ on a semiring $R$, we can define
  operations $+$ and $\cdot$ on its set of equivalence classes
  $R\modsim=\{[x]\mid x\in R\}$ by $[x]+[y]:=[x+y]$ and $[x]\cdot
  [y]:=[xy]$, for $x,y\in R$, turning $(R\modsim,+,\cdot)$ into a
  semiring, called the \emph{quotient semiring}.
\end{rem}

Note that if $R$ is a ring, there is a one-to-one correspondence
between congruences and ideals by identifying a congruence with its
$0$-class. Hence a ring is congruence-simple if and only if it is
simple in the sense that there are only trivial ideals.

By a semiring homomorphism we mean a map $f:R\to S$ between semirings
$R$ and $S$ which preserves the semiring operations and the zero
element.  Note that any homomorphism $f:R\to S$ gives rise to a
congruence $\sim$ on $R$ by defining $x\sim y$ if and only if
$f(x)=f(y)$, for $x,y\in R$. On the other hand, for any congruence
$\sim$ on $R$ we have the natural homomorphism $R\to R\modsim$. This
easily proves the following

\begin{rem}
  A semiring $R$ is congruence-simple if and only if any nonzero
  homomorphism $f:R\to S$ into a semiring $S$ is injective.
\end{rem}

The following example of a semiring turns out to be important.

\begin{exmp}
  Let $(M,+)$ be a commutative monoid. We call a map $f:M\to M$ an
  endomorphism if it preserves the monoid operation and the neutral
  element. On the set $\End(M)$ of all endomorphisms of $M$ we get
  operations $+$ and $\circ$ by defining $f+g$ as pointwise addition
  and and $f\circ g$ as composition of maps, for $f,g\in\End(M)$.
  
  It is straight-forward to verify that $(\End(M),+,\circ)$ is a
  semiring with a one, which will be called \emph{endomorphism
    semiring}.
\end{exmp}

The classification result uses subsemirings of some endomorphism
semirings, which are rich or lie dense in the sense that they contain
at least certain elementary endomorphisms.
  
\begin{defn}\label{defeab}
  Let $M$ be an idempotent commutative monoid.  A subsemiring
  $S\subseteq\End(M)$ is called \emph{dense} if it contains for all
  $a,b\in M$ the endomorphism $e_{a,b}\in\End(M)$, defined by
 \[e_{a,b}(x):=\begin{cases} 0&\text{if }x+a=a \\
    b&\text{otherwise}\end{cases}\qquad (x\in M).\]
\end{defn}

Now we can state the main result.
  
\begin{thm}\label{mr}
  Let $R$ be a finite semiring which is not a ring. Then the
  following are equivalent:
  \begin{enumerate}
  \item $R$ is congruence-simple.
  \item $|R|\leq 2$ or $R$ is isomorphic to a dense subsemiring
    $S\subseteq\End(M)$, where $(M,+)$ is a finite idempotent
    commutative monoid.
  \end{enumerate}
\end{thm}

Note that the classification of finite simple rings is a classical
subject in algebra. By the Wedderburn--Artin Theorem
(see~\cite{he68}), a finite ring $R$ with nontrivial multiplication
is simple if and only if $R$ is isomorphic to the endomorphism ring
$\Mat_{n\times n}(\F)$ of a finite-dimensional vector space $\F^n$
over a finite field $\F$.

\begin{rem}\label{remtwo}
  There are two proper semirings of order $2$, namely the semirings
  $R_{2,a}$, $R_{2,b}$ given by

  \[R_{2,a}:\ \text{
    \begin{tabular}{c|cc}
      $+$ & $0$ & $1$\\\hline
      $0$ & $0$ & $1$\\
      $1$ & $1$ & $1$
    \end{tabular}\quad
    \ \begin{tabular}{c|cc}
      $\cdot$ & $0$ & $1$\\\hline
      $0$ & $0$ & $0$\\
      $1$ & $0$ & $0$
    \end{tabular}}\qquad
  R_{2,b}:\ \text{
    \begin{tabular}{c|cc}
      $+$ & $0$ & $1$\\\hline
      $0$ & $0$ & $1$\\
      $1$ & $1$ & $1$
    \end{tabular}\quad
    \begin{tabular}{c|cc}
      $\cdot$ & $0$ & $1$\\\hline
      $0$ & $0$ & $0$\\
      $1$ & $0$ & $1$
    \end{tabular}}.\]
  
  $R_{2,b}$ is called the \emph{Boolean semiring} and can also be seen
  as the endomorphism semiring $\End(L_2)$ for $(L_2,+)
  =(\{0,1\},\max)$.  Trivially, $R_{2,a}$ and $R_{2,b}$ are
  congruence-simple.
\end{rem}

The proof of the main result is given in Sections~\ref{easy}
and~\ref{hard} of this paper.  In Section~\ref{easy} we show the
direction $(2)\Rightarrow (1)$, whereas in Section~\ref{hard} we
establish the direction $(1)\Rightarrow(2)$ with the help of
irreducible semimodules. Finally, we take a closer look at the
dense subsemirings of the endomorphism semirings in Section~\ref{mcsr}.

  
\section{Endomorphism semirings}\label{easy}

In this section we shall prove the direction $(2)\Rightarrow (1)$ of
Theorem~\ref{mr}. We begin with a remark on idempotent commutative
monoids and (semi-)lattices (see e.g.~\cite[sec. I.5 and II.2]{bi67}).

\begin{rem}\label{remlatt}
  Let $(M,+)$ be an idempotent commutative monoid. By defining $x\leq
  y$ if and only if $x+y=y$ for $x,y\in M$, we get a partial order
  relation $\leq$ on $M$, where $0\leq x$ for any $x\in M$. Also, for
  all $x,y\in M$ there exists a supremum $x\vee y=x+y$, so that
  $(M,\vee)$ is a join-semilattice.
    
  If in addition $M$ is finite, for all $x,y\in M$ there exists an
  infimum $x\wedge y=\sum_{z\leq x,\,z\leq y}z$, so that
  $(M,\vee,\wedge)$ is even a lattice.
  
  Now if $M$ is viewed as a lattice, the elements $f\in\End(M)$ are
  maps $f:M\to M$ satisfying $f(0)=0$ and $f(x\vee y)=f(x)\vee f(y)$
  for all $x,y\in M$. In particular, $f$ is order-preserving. Note
  however that $f(x\wedge y)=f(x)\wedge f(y)$ is not generally true,
  i.e. $f$ may not be a lattice endomorphism.
\end{rem}

Now we state a lemma on the maps $e_{a,b}$ of Definition~\ref{defeab}. Note
that by Remark~\ref{remlatt} we have \[e_{a,b}(x)=
\begin{cases} 
  0&\text{if }x\leq a\\ 
  b&\text{otherwise}
\end{cases}
\qquad(a,b,x\in M).\]

\begin{lem}\label{lemeab}    
  For $a,b\in M$, we have $e_{a,b}\in\End(M)$. Also, for $f\in\End(M)$
  and $a,b,c,d\in M$, we have $f\circ e_{a,b}=e_{a,f(b)}$ and
  \[e_{c,d}\circ f\circ e_{a,b}=
  \begin{cases}
    0&\text{if }f(b)\leq c,\\
    e_{a,d}&\text{otherwise}. 
  \end{cases}\] If $(M,+)$ has an
  absorbing element $\infty\in M$, then $e_{0,\infty}$ is absorbing
  for $(\End(M),+)$.
\end{lem}
  
\begin{proof}
  Note that for all $x,y\in M$, we have $x\vee y\leq a$ if and only if
  $x\leq a$ and $y\leq a$. It follows that $e_{a,b}(x\vee y)=0$ if and
  only if $e_{a,b}(x)=0$ and $e_{a,b}(y)=0$, that is if and only if
  $e_{a,b}(x)\vee e_{a,b}(y)=0$. Thus $e_{a,b}\in\End(M)$.
  
  Now if $f\in\End(M)$ and $a,b\in M$ one easily verifies $f\circ
  e_{a,b}=e_{a,f(b)}$. Applying this formula twice yields
  \[e_{c,d}\circ f\circ e_{a,b}=e_{c,d}\circ e_{a,f(b)}=e_{a,e_{c,d}(f(b))}=
  \begin{cases} 
    0&\text{if }f(b)\leq
    c,\\ e_{a,d}&\text{otherwise}.
  \end{cases}\] Finally, for any
  $h\in\End(M)$ and $x\in M\setminus\{0\}$ we have
  $(h+e_{0,\infty})(x)=h(x)+\infty=\infty$, so that
  $h+e_{0,\infty}=e_{0,\infty}$.
\end{proof}
  
\begin{prop}\label{propeab}
  Let $(M,+)$ be an idempotent commutative monoid with an absorbing
  element. Then any dense subsemiring $R\subseteq\End(M)$ is
  congruence-simple. In particular, $\End(M)$ itself is congruence-simple.
\end{prop}

Note that any \emph{finite} idempotent commutative monoid $M$ has an
absorbing element, namely $\infty:=\sum_{x\in M}x$.

\begin{proof}
  Let $\sim\;\subseteq R\times R$ be a semiring congruence relation.
  Suppose that $\sim\;\neq\id_R$, so that there exists $f,g\in R$ with
  $f\neq g$, but $f\sim g$. There is $b\in M$ with $f(b)\neq g(b)$,
  and without loss of generality, we may assume $f(b)\not\leq
  c:=g(b)$.
  
  For all $a,d\in M$ we have $e_{a,b}\in R$ and $e_{c,d}\in R$. Hence,
  since $\sim$ is a congruence, \[e_{c,d}\circ f\circ e_{a,b}\sim
  e_{c,d}\circ g\circ e_{a,b},\] so that $e_{a,d}\sim 0$, by
  Lemma~\ref{lemeab}.
  
  In particular $e_{0,\infty}\sim 0$, where $\infty\in M$ is the
  absorbing element. It follows that \[e_{0,\infty}=h+e_{0,\infty}\sim
  h+0=h\] for any $h\in R$, since $\sim$ is a congruence. Therefore
  $\sim\;=R\times R$, so that $R$ has no nontrivial congruence
  relations.
\end{proof}


\section{Finite congruence-simple semirings}\label{hard}

In this section we prove that any finite congruence-simple semiring
which is not a ring is of the form described in Theorem~\ref{mr}. We
start with a result established and proven by Monico in a more general
setting~\cite{mo04} and give a simplified proof for our case.

\begin{prop}\label{chris}
  Let $R$ be a congruence-simple semiring which is not a ring. Then
  the addition $(R,+)$ is idempotent.
\end{prop}

\begin{proof}
  For $x\in R$ and $n\in\N_0:=\{0,1,2,3,\dots\}$ let us write
  $nx:=x+\dots+x$, summing $x$ $n$-times. Also let $R+x:=\{y+x\mid
  y\in R\}$. Now, for $x,y\in R$ define \[x\sim y\quad
  :\Leftrightarrow\quad \exists\,m,n\in\N_0:\ mx\in R+y,\ ny\in R+x.\]
  Then it is easily verified that $\sim$ is a congruence relation.
  
  By congruence-simplicity it follows that $\sim\;=\id_R$ or
  $\sim\;=R\times R$. In the first case, since $x\sim x+x$, we deduce
  that $(R,+)$ is idempotent. In the second case, for all $x\in R$, we
  have $x\sim 0$, so that $0\in R+x$. This shows that $(R,+)$ is a
  group and thus $R$ is a ring.
\end{proof}

\begin{rem}\label{remchris}
  A congruence-simple semiring $R$ with idempotent addition and
  trivial multiplication $R\,R=\{0\}$ has order $\leq 2$.  Indeed,
  since $(R,+)$ is idempotent, $x+y=0$ implies $x=y=0$ for $x,y\in R$,
  so the equivalence relation $\sim$ on $R$ with classes $\{0\}$ and
  $R\setminus\{0\}$ is a congruence. Thus $\sim\;=\id_R$ and hence
  $|R|\leq 2$.
\end{rem}

\subsection{Semimodules}

The concept of semimodules over semirings is well-known
(see~\cite{go99}). For the proof of the classification result, we show
that any finite congruence-simple semiring admits a semimodule which
is irreducible in a strong sense and then we derive consequences from
it.

To fix some notations, let $R$ be a semiring.
  
\begin{defn}
  A (left) \emph{semimodule} $M$ over $R$ is a commutative monoid
  $(M,+)$ with neutral element $0\in M$, together with an
  $R$-multiplication $R\times M\to M$, $(r,x)\mapsto r\cdot x=rx$,
  such that, for all $r,s\in R$ and $x,y\in M$, we have $\ r(sx)=(rs)x,\ 0x=0,\ r0=0,\ $ and $\ (r+s)x=rx+sx,\ r(x+y)=rx+ry$.
\end{defn}
  
\begin{rem}
  If $(M,+)$ is a commutative monoid, any representation i.e. semiring
  homomorphism
  \[T:R\to\End(M),\quad r\mapsto T_r\] turns $M$ into a semimodule
  by defining $rx:=T_r(x)$, for $x\in R$ and $x\in M$.
    
  On the other hand, let $M$ be any semimodule over $R$. For $r\in R$,
  the map $x\mapsto rx$ defines an endomorphism $T_r$ of $M$, and the
  map $T:R\to\End(M)$, $r\mapsto T_r$ is a representation.
\end{rem}

\begin{defn}
  Let $M$ be a semimodule over $R$. A \emph{subsemimodule} $N\subseteq
  M$ is a submonoid of $(M,+)$ such that $R\cdot N\subseteq N$.  An
  equivalence relation $\sim$ on $M$ is called \emph{congruence} if
  \[x\sim y\quad\text{implies}\quad a+x\sim a+y,\ rx\sim
  ry,\quad\text{for all }x,y,a\in M\text{ and }r\in R.\]
\end{defn}

\begin{rem}
  Note that any subsemimodule $N\subseteq M$ itself is a semimodule
  over $R$. Also, given a congruence $\sim$ on $M$, we can define an
  addition and an $R$-multiplication on its set of equivalence classes
  $M\modsim=\{[x]\mid x\in M\}$ by \[[x]+[y]:=[x+y],\ 
  r[x]:=[rx],\quad\text{for all }x,y\in M,\ r\in R,\] turning
  $M\modsim$ into a semimodule over $R$, called the \emph{quotient
    semimodule}.
\end{rem}

If $M$ is a semimodule over $R$, let us call the subsemimodules
$\{0\}$ and $M$ and also the quotient semimodules $M/\id_M\cong M$ and
$M/(M\times M)\cong\{0\}$ the trivial ones.

\begin{defn}
  A semimodule $M$ over $R$ satisfying $R\,M\neq\{0\}$ is called 
  \begin{itemize}
  \item \emph{sub-irreducible} if it has only trivial subsemimodules,
  \item \emph{quotient-irreducible} if it has only trivial quotient
    semimodules,
  \item \emph{irreducible} if it is both sub-irreducible and
    quotient-irreducible.
  \end{itemize}
\end{defn}

Some authors refer to sub-irreducible and quotient-irreducible
semimodules as minimal and simple semimodules, respectively.

By a semimodule homomorphism we mean a map $f:M\to N$ between
semimodules over $R$ which preserves the semimodule operations as well
as the zero element. In this case, $f(M)$ is a subsemimodule of $N$,
and the relation $x\sim y$ if and only if $f(x)=f(y)$, for $x,y\in M$,
is a congruence on $M$. On the other hand, for any subsemimodule
$N_0\subseteq N$ and any quotient semimodule $M\modsim$ there are
natural homomorphisms $i:N_0\to N$ and $p:M\to M\modsim$.  This
constitutes the following

\begin{rem}\label{rem1}
  Let $M$ be a semimodule over $R$ such that $R\,M\neq\{0\}$. Then
  \begin{itemize}
  \item $M$ is sub-irreducible if and only if any nonzero
    homomorphism $f:N\to M$ from a semimodule $N$ is surjective,
  \item $M$ is quotient-irreducible if and only if any nonzero
    homomorphism $f:M\to N$ into a semimodule $N$ is injective.
  \end{itemize}
\end{rem}

\begin{rem}\label{rem2}
  To illustrate the use of irreducible semimodules we give a version
  of Schur's Lemma (see~\cite{he68}): Let $M$ be an irreducible
  semimodule over $R$ with representation $T:R\to\End(M)$, $r\mapsto
  T_r$. Then the commuting semiring \[C(M):=\{f\in\End(M)\mid f\circ
  T_r=T_r\circ f\text{ for all }r\in R\}\] is a semifield, i.e. any
  nonzero element is invertible.  Indeed, if $f\in
  C(M)\setminus\{0\}$, then $f:M\to M$ is a nonzero semimodule
  homomorphism, which by Remark~\ref{rem1} must be injective and
  surjective. It then easily follows that the inverse $f^{-1}$ lies in
  $C(M)$.
  
  In particular, if $(M,+)$ is finite and idempotent, then $C(M)$ is a
  finite proper semifield. It follows (see~\cite[sec. I.5]{he93a})
  that $C(M)$ has order $\leq 2$, so that $C(M)=\{0,\id_M\}$ is
  trivial.  If the representation $R\to\End(M)$ is faithful i.e.
  injective (this holds for example if $R$ is congruence-simple and
  $R\,M\neq\{0\}$), it follows that $R$ has trivial center, since
  \[\{x\in R\mid xr=rx\text{ for all }r\in
  R\}=T^{-1}(C(M))=\{0,1\}\cap R.\]
\end{rem}

\subsection{Existence of irreducible semimodules}

\begin{prop}\label{irrmod}
  Any finite congruence-simple semiring $R$ with $R\,R\neq\{0\}$
  admits a finite irreducible semimodule.
\end{prop}

To prove this result we begin with two lemmas that guarantee the
property $R\,M\neq\{0\}$ for certain semimodules $M$ over $R$. By a
nontotal semimodule congruence on $M$ we mean a congruence $\sim\;\neq
M\times M$, so that $M\modsim\;\neq\{0\}$. 

\begin{lem}\label{lemirr1}
  Let $R$ be a congruence-simple semiring with $R\,R\neq\{0\}$,
  considered as a semimodule over itself, and let $\sim$ be a nontotal
  semimodule congruence on $R$. Then, for the quotient semimodule
  $M:=R\modsim$ we have $R\,M\neq\{0\}$.
\end{lem}

\begin{proof}
  Since $\sim$ is a semimodule congruence, $r\sim s$ implies $x+r\sim
  x+s$ and $xr\sim xs$, for any $r,s,x\in R$. Now suppose
  $R\,M=\{0\}$.  Then for any $r,x\in R$ we have $[rx]=r[x]=0$, so
  that $rx\sim 0$.  Hence $r\sim s$ implies also $rx\sim sx$, for any
  $r,s,x\in R$, so that $\sim$ is even a semiring congruence. Since
  $\sim$ is nontotal, we must have $\sim\;=\id_R$ by
  congruence-simplicity. Hence $M=R$ and $R\,R=\{0\}$, which
  contradicts our assumption.
\end{proof}

\begin{lem}\label{lemirr2}
  Let $M$ be a semimodule over $R$ such that $R\,M\neq\{0\}$.
  \begin{enumerate}
  \item If $M$ is sub-irreducible, then $R\,P\neq\{0\}$ for all its
    nonzero quotient semimodules $P=M\modsim$.
  \item If $M$ is quotient-irreducible, then $R\,N\neq\{0\}$ for all
    its nonzero subsemimodules $N\subseteq M$.
  \end{enumerate}
\end{lem}
  
\begin{proof}
  (1) Let $M$ have only trivial subsemimodules. Since $R\,M\subseteq
  M$ is a subsemimodule, we must have $R\,M=M$. Now let $P=M\modsim$
  be a quotient subsemimodule with $R\,P=\{0\}$. Then we have
  $M=R\,M\subseteq[0]_\sim$, and therefore $M\modsim=\{0\}$.
    
  (2) Let $A:=\{x\in M\mid R\,x=\{0\}\}\subseteq M$ be the annulator
  of $R$ in $M$. Then it is easy to check that $A$ is a semimodule of
  $M$ with the additional property that $x\in A$ and $x+y\in A$
  implies $y\in A$. Also it is straightforward to check that defining
  \[x\sim y\quad:\Leftrightarrow\quad\exists\,a,b\in A:x+a=y+b\] for
  $x,y\in M$ gives a congruence $\sim$ on $M$ such that its zero-class
  $\{x\in M\mid x\sim 0\}$ equals $A$. Finally note that $A\neq M$ by
  assumption.
  
  Now if $M$ has only trivial quotient semimodules, the relation
  $\sim$ above must equal $\id_M$, and hence $A=\{0\}$. It follows
  that any subsemimodule $N\subseteq M$ with $R\,N=0$ must be zero.
\end{proof}

\begin{proof}[Proof of Proposition~\ref{irrmod}]
  We recursively define a sequence $M_0,M_1,\dots,M_n$ of finite
  semimodules over $R$ of decreasing sizes such that
  \begin{itemize}
  \item for all $i=0,\dots,n$ we have $R\,M_i\neq\{0\}$,
  \item for all $i=1,\dots,n$ we have $M_i$ is sub-irreducible or
    quotient-irreducible,
  \item $M_n$ is irreducible.
  \end{itemize}
  
  We start with $M_0:=R$, so that $R\,M_0=R\,R\neq\{0\}$.
  
  Now let $\sim$ be a maximal nontotal semimodule congruence on $R$
  (probably $\sim\;=\id_R$) and let $M_1:=R\modsim$. Since $\sim$ is
  nontotal we have $R\,M_1\neq\{0\}$ by Lemma~\ref{lemirr1}. By
  maximality of $\sim$ it follows that $M_1$ is quotient-irreducible.
  
  Suppose that $M_i$ has been defined for some $i\geq 1$, so that
  $R\,M_i\neq\{0\}$ and $M_i$ is sub-irreducible or
  quotient-irreducible. If $M_i$ is even irreducible we set $n=i$ and
  stop.
  
  Otherwise suppose that $M_i$ is quotient-irreducible but has
  nontrivial subsemimodules. Take a minimal nonzero semimodule
  $M_{i+1}\subseteq M_i$. Then $R\,M_{i+1}\neq\{0\}$ by
  Lemma~\ref{lemirr2}, (2), and furthermore $M_{i+1}$ is
  sub-irreducible.  Now consider the case where $M_i$ is
  sub-irreducible but has nontrivial congruences. By taking a maximal
  nontotal congruence $\sim$ and letting $M_{i+1}:=M_i\modsim$, we
  have $R\,M_{i+1}\neq\{0\}$ by Lemma~\ref{lemirr2}, (1), and
  furthermore $M_{i+1}$ is quotient-irreducible.
  
  The sequence has been constructed. Since $R$ is finite and the
  cardinalities of $M_1,M_2,\dots$ are strictly decreasing the
  sequence must terminate by an irreducible semimodule $M_n$ over $R$.
\end{proof}

\subsection{A density result}

Let $R$ be a congruence-simple semiring and $M$ be a semimodule over
$R$ with $R\,M\neq\{0\}$. Then the representation $R\to\End(M)$ is
nonzero and hence must be injective, so that $R$ can be seen as a
subsemiring of $\End(M)$. If $M$ is irreducible the question of the
`density' of $R$ in $\End(M)$ arises. We have already seen in
Remark~\ref{rem2} that the commuting semiring of $R$ in $\End(M)$ is
trivial if $(M,+)$ is idempotent. Now we show another density result:

\begin{prop}\label{propfinish}
  Let $R$ be a finite congruence-simple semiring with idempotent
  addition and let $M$ be a finite irreducible semimodule over $R$.
  Then $(M,+)$ is idempotent, and for all $a,b\in M$ there exists
  $r\in R$ such that
  \[rx=\begin{cases} 0&\text{if }x+a=a\\ 
    b&\text{otherwise}\end{cases}\qquad(x\in M).\]
  Thus $R$, seen as a subsemiring of $\End(M)$, is dense (see
  Definition~\ref{defeab}).
\end{prop}

\begin{proof}
  First note that $(M,+)$ is idempotent: By irreducibility, the
  subsemimodule $R\,M$ of $M$ is nonzero, hence $R\,M=M$. So, any
  $x\in M$ can be written as $x=ry$ with $r\in R$ and $y\in M$. It
  follows $x+x=ry+ry=(r+r)y=ry=x$, since $(R,+)$ is idempotent, so
  that $(M,+)$ is idempotent. Recall from Remark~\ref{remlatt} that
  now on $M$ there is an order relation $\leq$ defined by $x\leq y$ if
  and only if $x+y=y$, for $x,y\in M$.
   
  For $x\in M$ define $I_x:=\{r\in R\mid rx=0\}$, which is a
  subsemimodule of $R$. We have $I_{x+y}=I_x\cap I_y$ for $x,y\in M$,
  since $rx+ry=0$ implies $rx=ry=0$ for $r\in R$, because $(M,+)$ is
  idempotent. Now we claim that defining \[x\sim y\quad
  :\Leftrightarrow\quad I_x=I_y\qquad (x,y\in M)\] gives a semimodule
  congruence on $M$: Indeed, if $x\sim y$ and $z\in M$, we have
  $I_{z+x}=I_z\cap I_x=I_z\cap I_y=I_{z+y}$, so that $z+x\sim z+y$.
  Also for $r,s\in R$ we have $r(sx)=(rs)x=0$ if and only if
  $(rs)y=r(sy)=0$, so that $I_{sx}=I_{sy}$ i.e. $sx\sim sy$.
  
  Assume that $\sim\;=M\times M$. Then $I_x=I_0=R$ for all $x\in M$,
  so that $R\,M=\{0\}$, which cannot hold. Since $M$ is
  quotient-irreducible it follows that $\sim\;=\id_M$. We conclude
  that $x\leq y$ is equivalent to $I_y\subseteq I_x$, for $x,y\in M$,
  since $x+y=y$ if and only if $I_x\cap I_y=I_{x+y}=I_y$.
  
  Now let $a\in M$ be fixed. If $a=\infty$, the absorbing element in
  $(M,+)$, the assertion trivially holds with $r=0$. So assume
  $a\neq\infty$. For any $x\in M$ with $x\not\leq a$ we have shown
  before that $I_a\not\subseteq I_x$, so the semimodule homomorphism
  $I_a\to M$, $r\mapsto rx$ is nonzero. Since $M$ is sub-irreducible,
  it must be surjective, so in particular there exists $r_x\in I_a$
  such that $r_x\,x=\infty$. Letting $s:=\sum_{x\not\leq a}r_x\in
  I_a\subseteq R$, for $x\in M$ we have \[sx=\begin{cases}0&\text{if
      }x\leq a,\text{ since then }sx=sx+sa=sa=0,\\\infty&\text{if
      }x\not\leq a,\text{ since then }sx\geq
    r_x\,x=\infty,\end{cases}\] so we have shown the assertion for
  $b=\infty$.
  
  Consider now the subsemimodule $N:=\{r\infty\mid r\in R\}$ of $M$.
  We have $\infty=s\infty\in N$, so that $N\neq\{0\}$. By
  sub-irreducibility of $M$ it follows $N=M$, so for any $b\in M$
  there exists $r\in R$ with $r\infty=b$. Then for $x\in M$ we have
  $(rs)x=0$ if $x\leq a$, and $(rs)x=b$ otherwise, which completes the
  proof.
\end{proof}

Now we complete the proof of the Theorem~\ref{mr} by showing the
direction $(1)\Rightarrow (2)$. Let $R$ be a finite congruence-simple
semiring which is not a ring and suppose $|R|>2$. Then $(R,+)$ is
idempotent by Proposition~\ref{chris} and $R\,R\neq\{0\}$ by
Remark~\ref{remchris}.  Afterwards, Proposition~\ref{irrmod} guarantees
the existence of a finite irreducible semimodule $M$ over $R$, so that
$R$ is isomorphic to a subsemiring $S$ of $\End(M)$. Finally, by
Proposition~\ref{propfinish} we have that $S$ is a dense subsemiring
of $\End(M)$.


\section{The family of dense endomorphism subsemirings}\label{mcsr}

\begin{defn}
Let $M$ be an idempotent commutative monoid. We define $\mcSR(M)$ to be the collection of all dense subsemirings $R\subseteq\End(M)$.
\end{defn}

In this section we take a closer look at the families $\mcSR(M)$.
First we address the question of isomorphy and anti-isomorphy of these
semirings. Then we give a criterion when the family $\mcSR(M)$ is
trivial. Finally we list the dense endomorphism subsemirings having
smallest order.

In this section, let $M,M_1$ and $M_2$ be always idempotent
commutative monoids having an absorbing element.

\subsection{Isomorphy}

\begin{prop}\label{propiso}
  Let $R_1\in\mcSR(M_1)$ and $R_2\in\mcSR(M_2)$ be isomorphic
  semirings. Then also the monoids $M_1$ and $M_2$ are isomorphic.
\end{prop}

We prove a lemma first.  Recall from Lemma~\ref{lemeab} that if
$\infty\in M$ is the absorbing element, then $e_{0,\infty}$ is an
absorbing element in $(R,+)$ for any semiring $R\in\mcSR(M)$.

\begin{lem}\label{lemiso}
  Let $R\in\mcSR(M)$ and let $z\in R$ be the absorbing element in
  $(R,+)$. Then the map \[\theta:M\to Rz,\quad b\mapsto e_{0,b}\]
  defines an isomorphism between $(M,+)$ and the submonoid $Rz$ of
  $(R,+)$.
\end{lem}

\begin{proof}
  Note that $f\circ e_{0,\infty}=e_{0,f(\infty)}$ for all $f\in R$, so
  in particular $e_{0,b}\circ e_{0,\infty}=e_{0,b}$ for all $b\in M$.
  This shows $Rz=Re_{0,\infty}=\{e_{0,b}\mid b\in M\}$, so $\theta$ is
  well-defined and surjective. That $\theta$ is injective and a
  homomorphism is clear.
\end{proof}

\begin{pf*}{Proof of Proposition~\ref{propiso}}
  Suppose there is a semiring isomorphism $\phi:R_1\to R_2$. For
  $i=1,2$, let $z_i\in R_i$ be the absorbing element in $(R_i,+)$. We
  then have $\phi(z_1)=z_2$ and thus $\phi(R_1z_1)=R_2z_2$. The
  restriction $\phi'=\phi|_{R_1z_1}:R_1z_1\to R_2z_2$ of $\phi$ is
  therefore an isomorphism between the submonoids $R_1z_1$ and
  $R_2z_2$ of $(R_1,+)$ and $(R_2,+)$, respectively. Now for $i=1,2$,
  let $\theta_i:M_i\to R_iz_i$ be the isomorphism defined in
  Lemma~\ref{lemiso}. Then we can construct an isomorphism
  \[\theta_2^{-1}\circ\phi'\circ\theta_1:M_1\to M_2\] between the
  monoids $(M_1,+)$ and $(M_2,+)$.
\end{pf*}

Next we identify anti-isomorphic pairs of congruence-simple semirings.

\begin{rem}
  Let $M$ be finite with corresponding lattice $(M,\vee,\wedge)$, so
  that $(M,+)=(M,\vee)$. Then also $(M,\wedge)$ is a finite idempotent
  commutative monoid, which we denote by $\tilde{M}$. Its
  corresponding lattice is the \emph{dual lattice} of $M$, obtained by
  reversing the order $(M,\leq)$. 
\end{rem}
  
Let $(L_2,\vee)=(\{0,1\},\max)$ and let $M^*=\Hom(M,L_2)$ be the set of
all monoid homomorphisms $M\to L_2$. Defining addition pointwise, $M^*$
becomes a finite idempotent commutative monoid.

\begin{lem}\label{lemdual}
  The monoid $M^*$ is isomorphic to $\tilde{M}$. In fact, the map
  \[M\to M^*,\quad a\mapsto e_a,\ \text{where }e_a(x)=
  \begin{cases}
    0&\text{if }x\leq a,\\
    1&\text{otherwise,}
  \end{cases}\] is a bijection such that
  $e_{a\wedge b}=e_a\vee e_b$ for all $a,b\in M$. 
\end{lem}

\begin{proof}
  This is rephrasing the well-known result in lattice theory that any
  finite lattice is isomorphic to its lattice of ideals
  (see~\cite[sec. II.3]{bi67}).
\end{proof}

\begin{prop}\label{propdual}
  Let $M$ be finite. The semirings $\End(M)$ and $\End(\tilde{M})$ are
  anti-isomorphic.
\end{prop}

\begin{proof}
  By Lemma~\ref{lemdual} we may assume $\tilde{M}=M^*$. Define a map
  \[\End(M)\to\End(M^*),\quad f\mapsto f^*,\ \text{where
    }f^*(\phi):=\phi\circ f\text{ for }\phi\in M^*.\] It is easy to
  see that this map is well-defined and that the following algebraic
  properties hold for $f,g\in\End(M)$: \[(f+g)^*=f^*+g^*,\quad
  0^*=0,\quad (f\circ g)^*=g^*\circ f^*.\] To prove injectivity,
  suppose we have $f,g\in\End(M)$ with $f^*=g^*$. With $e_a$ as
  defined in Lemma~\ref{lemdual} it follows $e_a(f(x))=e_a(g(x))$ for
  all $a,x\in M$, so that $f(x)\leq a$ if and only if $g(x)\leq a$.
  For all $x\in M$ it follows $f(x)=g(x)$, hence $f=g$.
  
  From injectivity it follows in particular
  $|\End(M)|\leq|\End(\tilde{M})|$. We can apply this result to
  $\tilde{M}$ to yield $|\End(\tilde{M})|\leq|\End(M)|$. Thus
  $|\End(M)|=|\End(\tilde{M})|$ and the map is also surjective.
\end{proof}

\begin{cor}\label{cordual}
  Let $M$ be finite and suppose $M$ as a lattice is isomorphic to its
  dual lattice. Then the semiring $\End(M)$ is anti-isomorphic to
  itself.
\end{cor}

\begin{cor}
  Let $M_1$ and $M_2$ be finite and let $R_1\in\mcSR(M_1)$ and
  $R_2\in\mcSR(M_2)$ be anti-isomorphic semirings. Then the monoids
  $M_1$ and $\tilde{M_2}$ are isomorphic.
\end{cor}

\begin{proof}
  By Proposition~\ref{propdual}, $\End(M_2)$ is anti-isomorphic to
  $\End(\tilde{M_2})$, and thus $R_1$ is isomorphic to some
  $R_2'\in\mcSR(\tilde{M_2})$. Now the result follows from
  Proposition~\ref{propiso}.
\end{proof}

\subsection{The case $|\mcSR(M)|=1$}

We now discuss under which circumstances the only dense subsemiring of
$\End(M)$ is $\End(M)$ itself.

\begin{prop}\label{propone}
  Let $M$ be finite. Then we have $\mcSR(M)=\{\End(M)\}$ if and only
  if the lattice $(M,\vee,\wedge)$ satisfies the following condition:
  \[\tag{D} \forall z\in M:\ z=\bigvee_{a,\,z\not\leq a}\ 
  \bigwedge_{x,\,x\not\leq a}x.\]
\end{prop}

\begin{proof}
  If $S$ is the subsemiring of $R:=\End(M)$ generated by the set
  $E:=\{e_{a,b}\mid a,b\in M\}$, then we have $\mcSR(M)=\{\End(M)\}$
  if and only if $S=R$. Note that since $E$ is closed under
  multiplication (see Lemma~\ref{lemeab}) $S$ consists of all finite
  sums of elements in $E$. Writing $1=\id_M\in R$ we show that
  \[\tag{$*$}S=R\quad\text{if and only if}\quad 1=\sum_{(a,b)\in
    X}e_{a,b}\text{ with }X:=\{(a,b)\in M^2\mid e_{a,b}\leq1\}.\]
  
  Indeed, suppose $S=R$, so we can express in particular $1$ as a sum
  of elements in $E$, say $1=\sum_ie_{a_i,b_i}$. Surely,
  $e_{a_i,b_i}\leq1$ and hence $(a_i,b_i)\in X$ for all $i$, so that
  \[1=\sum_ie_{a_i,b_i}\leq\sum_{(a,b)\in X}e_{a,b}\leq1\] and thus
  the right side of $(*)$ holds. On the other hand, supposing
  $1=\sum_{(a,b)\in X}e_{a,b}$ implies $1\in S$. Then for any $f\in R$
  we have \[f=f\circ1=\sum_{(a,b)\in X}f\circ e_{a,b}=\sum_{(a,b)\in
    X}e_{a,f(b)}\in S\] (see Lemma~\ref{lemeab}), so that $S=R$. This
  proves the equivalence $(*)$.
  
  Note next that $(a,b)\in X$ i.e. $e_{a,b}\leq1$ if and only if
  $b\leq x$ for all $x\not\leq a$ which is equivalent to
  $b\leq\bigwedge_{x,\,x\not\leq a}x$. This shows that
  \[\sum_{(a,b)\in X}e_{a,b}=\sum_{a\in M}e_{a,b_a}\text{ with
    }b_a:=\bigwedge_{x,\,x\not\leq a}x.\]
  
  Now for all $z\in M$ we have \[\sum_{(a,b)\in
    X}e_{a,b}(z)=\sum_{a\in M}e_{a,b_a}(z)=\bigvee_{a,\,z\not\leq
    a}b_a=\bigvee_{a,\,z\not\leq a}\ \bigwedge_{x,\,x\not\leq a}x,\]
  which together with $(*)$ concludes the proof.
\end{proof}

\begin{rem}
  The condition (D) given in proposition~\ref{propone} is fulfilled if
  and only if the lattice $M$ is distributive, or equivalently, $M$ is
  isomorphic to a ring of subsets (cf.~\cite[sec. III.3]{bi67}).
  
  Indeed, assume that $(M,\cup,\cap)$ is a ring of subsets, i.e. a
  sublattice of a power set lattice $(\mc{P}(\Omega),\cup,\cap)$. For
  $\omega\in\Omega$ let $A_{\omega}:=\bigcup_{X\in M,\omega\notin
    X}X\in M$. Then for $X\in M$ we have $X\subseteq A_{\omega}$ if
  and only if $\omega\notin X$. It follows
  \[Z\supseteq\bigcup_{A,\,Z\not\subseteq A}\
  \bigcap_{X,\,X\not\subseteq A}X
  \supseteq\bigcup_{\omega,\,Z\not\subseteq A_{\omega}}\
  \bigcap_{X,\,X\not\subseteq
    A_{\omega}}X=\bigcup_{\omega,\,\omega\in Z}\
  \bigcap_{X,\,\omega\in X}X\supseteq Z\] for all $Z\in M$, so $M$
  satisfies property (D).
  
  On the other hand, if we have a lattice $(M,\vee,\wedge)$ with
  condition (D), let $\Omega:=\{b_a\mid a\in M\}$ with
  $b_a:=\bigwedge_{x,\,x\not\leq a}x$. Consider the representation of
  $M$ given by \[\Phi:M\to\mc{P}(\Omega), \quad z\mapsto\{b_a\mid a\in
  M,\,z\not\leq a\}.\]
  
  We can see directly that $z_1\leq z_2$ implies
  $\Phi(z_1)\subseteq\Phi(z_2)$. On the other hand, with the help of
  (D) we conclude that $\Phi(z_1)\subseteq\Phi(z_2)$ implies
  $z_1=\bigvee_{a,\,z_1\not\leq a}b_a\leq \bigvee_{a,\,z_2\not\leq
    a}b_a=z_2$. It follows that $\Phi$ is a lattice monomorphism, so
  that $M$ is isomorphic to a sublattice of
  $(\mc{P}(\Omega),\cup,\cap)$.
\end{rem}

\subsection{Congruence-simple semirings of small order}

Table~\ref{tab:1} shows the smallest nontrivial idempotent commutative
monoids $M$ (up to isomorphy), represented by the Hasse-diagram of the
corresponding lattices, together with the semirings in the collection
$\mcSR(M)$. We write $R_m$ for a semiring with $m$ elements.

These, together with $R_{2,a}$ from Remark~\ref{remtwo}, are the
smallest congruence-simple semirings which are not rings. The smallest
such semiring not shown in Table~\ref{tab:1} has order $98$.

Note that $R_{50,a}$ and $R_{50,b}$ are anti-isomorphic to each other
by Proposition~\ref{propdual}, whereas the other semirings in
Table~\ref{tab:1} are self-anti-isomorphic by Corollary~\ref{cordual}.
Furthermore, all semirings in Table~\ref{tab:1} have a one-element, except
$R_{42}$ and $R_{44}$.

\begin{ack}
  I am very grateful to my advisor Joachim Rosenthal and also to
  G\'{e}rard Maze for their great encouragement and some fruitful
  discussions.
\end{ack}

\renewcommand{\arraystretch}{1.5}

\begin{table}
\caption{The smallest lattices together with the corresponding 
endomorphism semirings.}
\label{tab:1}\vspace{5mm}

\begin{tabular}{cc}
  
  \begin{minipage}{6cm}
    \begin{tabular}{|c|p{3cm}|}
      \hline
      $M$ & $\mcSR(M)$ \\
      \hline\hline
       \xymatrix @=0.4cm {
        \node \edge{d} \\
        \node
        } \vspace{-3mm} & $\{R_{2,b}\}$

      (the Boolean semiring) \\ & \\
      \hline\hline
      \xymatrix @=0.4cm {
        \node \edge{d} \\
        \node \edge{d} \\
        \node
        } \vspace{-3mm} & $\{R_6\}$ \\ & \\
      \hline\hline
      \xymatrix @=0.4cm {
        \node \edge{d} \\
        \node \edge{d} \\
        \node \edge{d} \\
        \node
        } \vspace{-2mm} & $\{R_{20}\}$ \\ & \\
      \hline
      \xymatrix @=0.4cm {
        & \node \\
        \node \edge{ur}\edge{dr} & & \node \edge{ul}\edge{dl} \\
          & \node
        } \vspace{-2mm} & $\{R_{16}\}$

      (the $2\times 2$-matrices over $R_2$) \\ & \\
      \hline
    \end{tabular}
  \end{minipage} &
  
  \begin{minipage}{7cm}
    \begin{tabular}{|c|p{4cm}|}
      \hline
      $M$ & $\mcSR(M)$ \\
      \hline\hline
      \xymatrix @=0.3cm {
        \node \edge{d} \\
        \node \edge{d} \\
        \node \edge{d} \\
        \node \edge{d} \\
        \node
        } \vspace{-3mm} & $\{R_{70}\}$ \\ & \\
      \hline
      \xymatrix @=0.3cm {
        & \node \\
        \node \edge{ur}\edge{dr} & & \node \edge{ul}\edge{dl} \\
        & \node \edge{d} \\
        & \node
        } \vspace{-2mm} & $\{R_{50,a}\}$ \\ & \\
      \hline
      \xymatrix @=0.3cm {
        & \node \edge{d} \\
        & \node \\
        \node \edge{ur}\edge{dr} & & \node \edge{ul}\edge{dl} \\
        & \node
        } \vspace{-2mm} & $\{R_{50,b}\}$ \\ & \\
      \hline
      \xymatrix @R=0.08cm @C=0.3cm {
        & \node \\ \\
        \node\edge{uur}\edge{dd} \\
        & & \node\edge{uuul}\edge{dddl} \\
        \node\edge{ddr} \\ \\
        & \node
        } \vspace{-2mm} & $\{R_{43}$, $R_{42}\}$ \\ & \\
      \hline
      \xymatrix @=0.4cm {
        & \node \\
        \node\edge{ur}\edge{dr} & \node\edge{u}\edge{d} &
        \node\edge{ul}\edge{dl} \\
        & \node
        } \vspace{-2mm} & $\{R_{50,c}$, $R_{47}$, $R_{46,a}$, 
      $R_{46,b}$, $R_{46,c}$, $R_{45}$, $R_{44}\}$

      (where $R_{46,a}$, $R_{46,b}$ and $R_{46,c}$ are isomorphic) \\ & \\
      \hline
    \end{tabular}
  \end{minipage}
  
\end{tabular}\vspace{1cm}
\end{table}

\def\cprime{$'$} \def\cprime{$'$}
\providecommand{\bysame}{\leavevmode\hbox to3em{\hrulefill}\thinspace}
\providecommand{\MR}{\relax\ifhmode\unskip\space\fi MR }
\providecommand{\MRhref}[2]{%
  \href{http://www.ams.org/mathscinet-getitem?mr=#1}{#2}
}
\providecommand{\href}[2]{#2}

\end{document}